\documentclass{amsart}
\usepackage{amssymb,amsmath,amscd,graphicx,amsthm}

\usepackage{color}

\newcommand{\prin}{\textit{prin}}
\newcommand{\univ}{\textit{univ}}

\newcommand{\A}{\mathcal A}

\newcommand{\Spec}{\operatorname{Spec}}

\renewcommand{\dagger}{\vee}
\newcommand{\New}{\operatorname{New}}

\newcommand{\ttt}{\mathfrak r}
\newcommand{\red}{\mathit{red}}
\def\C{{\mathbb C}}
\def\FF{{\mathcal F}}

\def\d{\mathbf d}

\newtheorem{proposition}{Proposition}
\newtheorem{thm}{Theorem}
\newtheorem{remark}{Remark}

\title{Associahedra as moment polytopes}
\author{Michael Gekhtman}
\address{Department of Mathematics, University of Notre Dame, Notre Dame, IN 46556 USA}
\email{mgekhtma@nd.edu}
\author{Hugh Thomas}
\address{LACIM, Université du Québec à Montréal, Montréal, QC, Canada}
\email{thomas.hugh\_r@uqam.ca}

\begin{document}
\maketitle

\begin{abstract} Generalized associahedra are a well-studied family of polytopes associated to a finite-type cluster algebra and choice of starting cluster. We show that the generalized associahedra constructed by Padrol, Palu, Pilaud, and Plamondon, building on ideas from Arkani-Hamed, Bai, He, and Yan, can be naturally viewed as moment polytopes for an open patch of the quotient of the cluster $\mathcal A$-variety with universal coefficients by its maximal natural torus action. We prove our result by showing that the construction of Padrol, Palu, Pilaud, and Plamondon can be understood on the basis of the way that moment polytopes behave under symplectic reduction.

  \end{abstract}

\section{Introduction}
Let $B$ be a skew-symmetrizable $n\times n$ matrix. The cluster algebra $A(B)$
is defined starting from $B$. We assume that the number of cluster variables of $A(B)$ is finite, and we list them as $z_1,\dots,z_v$ in some order.
The cluster variable $z_i$ has an associated $g$-vector which we denote
$g_i$. The $g$-vector fan $\Sigma(B)$ for $A(B)$ is the fan whose rays are
generated by the vectors
$g_i$, and whose maximal cones are given by the clusters of $A(B)$. 
By definition, a generalized associahedron for $A(B)$ is a polytope whose
outer normal fan is $\Sigma(B)$.

It seemed natural to us to wonder whether these polytopes arise naturally as moment polytopes. 
The goal of this paper is to answer this question
in the affirmative.  (At the end of this introduction, we shall explain some other routes to an affirmative answer in some cases.)

The approach which we shall take is as follows.
\cite{ABHY} constructs a (type $A_n$) associahedron in
a particular way, as the intersection of a non-negative orthant with an
affine subspace. (This construction was arrived at for purposes coming from the physics of scattering amplitudes, which will not concern us here.)
It turned out that this same construction could be applied more generally
\cite{BDMTY}.
In fact, it was shown in \cite{PPPP} that this construction
can
be extended to realize all associahedra of finite type cluster algebras.
In the construction of \cite{PPPP}, the associahedra for $A(B)$
are constructed
by intersecting $\mathbb R^v_{\geq 0}$ with
a certain $n$-dimensional affine subspace $\mathbb E$ (which depends on
$B$ and some freely chosen positive parameters).  
We will explain this construction in terms of symplectic geometry,
as an instance of symplectic reduction.

This interpretation is as follows.
We start with $X=(\mathbb C^*)^n \times \mathbb C^v$, which admits the obvious action by the torus $T=(\mathbb C^*)^{n+v}$.
An appropriate real symplectic form respects the torus action, and the  
image of $X$ under the corresponding moment map is $\mathbb R^n \times \mathbb R^v_{\geq 0}$.

There is a $v$-dimensional torus which acts on the cluster $\mathcal A$-variety for $B$ with universal coefficients. Viewing $X$ as a patch of that variety, the $v$-dimensional torus is a subtorus of $T$, which we denote $T_K$. 
Carrying out symplectic reduction with respect to $T_K(\mathbb R)$, we obtain a
variety whose moment map image, with respect to the residual torus $(T/T_K)(\mathbb R)$ is exactly
$\mathbb R^v_{\geq 0}\cap \mathbb E$.
The different choices of $\mathbb E$ correspond to different choices in the
symplectic reduction process.
The resulting
variety contains as a dense subset the cluster configuration
space of \cite{AHL}, which provides some motivation for considering this quotient.

The above approach provides a symplectic origin for all the associahedra constructed by \cite{PPPP} for finite type cluster algebras and arbitrary choice of initial seed.

\subsection{Other approaches to viewing associahedra as moment polytopes}

For $P$ a lattice polytope in $\mathbb R^n$, let $\Sigma_P$ be its outer normal fan. For any (rational) fan $\Sigma$, there is an associated toric variety
$X_\Sigma$. A lattice polytope $P$ equips $X_{\Sigma_P}$ with a map to projective
space, which can be used to define a moment map from $X_P$ to $\mathbb R^n$, whose image is $P$. (See \cite[Section 12.2]{CLS}.) Thus, for any lattice polytope whatsoever, it is possible to construct a variety whose moment map image is the polytope in question. Obviously, though, the fact that one can do this does not provide any explanation as to the origin of the lattice polytope.

A variant of this approach is applicable to associahedra in particular, as follows from 
Arkani-Hamed, He, and Lam \cite{AHL}.
They study the \emph{cluster configuration space} associated to any cluster
algebra $A(B)$ and choice of initial seed. They show that it embeds as a
dense set
inside the toric variety $X_{\Sigma(B)}$. Thus, given an associahedron $P$, it is
possible to proceed as in the previous paragraph to put a symplectic form on
the cluster configuration space
which will realize $P$ as its moment polytope, but, again, this is not
fully satisfactory as an explanation of where $P$ comes from.

A different approach which realizes (some) generalized associahedra is provided by Escobar \cite{Es}. Starting from $G$ a complex Lie group with Weyl group $W$, and a
word $q=s_{i_1}\dots s_{i_m}$ in the Coxeter generators, she defines the
\emph{brick manifold} as a certain submanifold of the Bott-Samelson
variety associated to $q$.
The brick manifold inherits a symplectic form from the Bott-Samelson
variety, such that the image of the moment map is the \emph{brick polytope} associated to $q$ 
in the sense of Pilaud--Stump \cite{PilaudStump}.
\cite{PilaudStump} showed that, for  $c$ a Coxeter element, if $q$ is chosen to be be equal to $c$
followed by the $c$-sorting word for $w_0$, then the
brick polytope is a realization of the associahedron with initial seed
corresponding to $c$. This procedure therefore realizes as moment polytopes the
finite type associahedra corresponding to acyclic seeds.

\section{Background}

\subsection{Cluster algebras}
Cluster algebras were discovered by Fomin and Zelevinsky \cite{CA1}. They
are now the subject of a vast literature. We will recall the small part of
the subject which we need.

An $n\times n$ integer matrix $B$ is said to be skew-symmetrizable if there is
a matrix $D$ with non-zero entries on the diagonal and zeros elsewhere, such
that $DB$ is skew-symmetric. From now on, we assume that $B$ is skew-symmetrizable. Note that if $B$ is skew-symmetrizable, so is the transpose of $B$, which we denote $B^\vee$.

Adding $m\geq 0$ more rows to $B$ if desired,
we obtain a matrix $\widetilde B$, from which we can define a cluster algebra
$A(\widetilde B)$.
This is a commutative algebra with a distinguished set of generators called
cluster variables.

For $1 \leq i \leq n$, we define a mutation operation $\mu_i$ on
$\widetilde B$ as
follows:

$$ \mu_i(\widetilde B)_{jk}= \left\{ \begin{array}{ll} -\widetilde B_{jk} & \textrm{if $j=i$ or $k=i$,}\\ \widetilde B_{jk}+[\widetilde B_{ji}]_+[\widetilde B_{ik}]_+ - [\widetilde B_{ji}]_-[\widetilde B_{ik}]_-&\textrm{otherwise.}\end{array}\right.$$
  Here $[x]_+=\max(x,0)$ and $[x]_-=\min(x,0)$.

  To construct the cluster algebra $A(\widetilde B)$, it is
  convenient to begin with an infinite $n$-regular tree $\mathbb T_n$
  with an edge labelling
  such that the edges incident to each vertex are labelled 1 to $n$. We pick a starting vertex, and
  associate to it the matrix $\widetilde B$, and the collection of variables
  $x_1,\dots,x_{n+m}$, inside the rational function field $\mathbb C(x_1,\dots,x_{n+m})$.

  We will now define a procedure which allows us to propagate this labelling
  to all vertices of $\mathbb T_n$. Suppose we have a vertex labelled
  $(\widetilde B', z_1,\dots,z_{n+m})$ which is incident to an edge
  labelled $i$, whose other end is not yet labelled. Define
  $$ z_i'= \frac{\prod_{\widetilde B'_{ki}>0} z_k^{\widetilde B'_{ki}} +  \prod_{\widetilde B'_{ki}<0} z_k^{-\widetilde B'_{ki}} }{z_i}$$

  We then label the node at the other end of the edge we are looking at by the
  matrix $\mu_i (\widetilde B')$ and the sequence of variables
  $z_1,z_2,\dots, z'_i,\dots,z_n$, where $z_i'$ has replaced $z_i$, but the
  other variables are unchanged.

  The cluster algebra $A(\widetilde B)$ is the $\mathbb C$-algebra generated by
  all the variables labelling all the nodes. (Actually we could use $\mathbb Z$ as the ground ring, but for our purposes it is natural to work over $\mathbb C$.)

  Note that the variables $x_{n+1},\dots,x_{n+m}$ will be appear at every node,
  because we only mutate at positions 1 to $n$. These are called
  ``frozen variables.''
  Depending on circumstances,
  it can be convenient to include their inverses in the cluster algebra as
  well. Sometimes it is convenient to distinguish them notationally by using
  a different letter for them. 

  We now consider some particular choices of $\widetilde B$ obtained by adding
  rows to an $n\times n$ matrix $B$.
  
One canonical way to add rows is to add a copy of the $n\times n$
identity matrix below $B$. The resulting matrix is called $B^\prin$.
The corresponding cluster algebra $A(B^\prin)$ is called the
cluster algebra for $B$ with principal coefficients, and it is conventional
to let the initial variables be $x_1,\dots,x_n,y_1,\dots,y_n$. This algebra
is
$\mathbb Z^n$-graded, with $x_i$ having grading
the standard basis vector $e_i$, while the degree of $y_i$ is the negative
of the vector obtained by reading the $i$-th column of $B$.
The degrees of the non-frozen cluster variables are
called the $g$-vectors associated to $B$. Distinct cluster variables
have distinct $g$-vectors.

Let us pause to note a few facts about the relationship between the cluster
algebra associated to $\widetilde B$ and that associated to $B$. Setting the
frozen variables to 1 will send a cluster variable for $\widetilde B$ to a
corresponding cluster variable for $B$. In principle, one could worry that
distinct cluster variables for $\widetilde B$ could be sent to the same
cluster variable for $B$, but this will not happen in the case of
interest to us, when $A(B)$ has only finitely many cluster variables. 
(Conjecturally it will never happen at all). Given a cluster variable for
$B$, we are therefore entitled to consider the corresponding variable in
$B^\prin$ (or $B^\univ$, to be introduced shortly) and therefore, in particular, there is a well-defined
$g$-vector for each cluster variable in $A(B)$.

\subsection{Universal cluster algebras}

Suppose that $A(B)$ has only a finite number $v$ of cluster variables. In this case there are the same number of cluster variables for $A(B^\vee)$. 

In this setting, we can define another matrix extending
$B$, which we denote $B^\univ$. It is obtained by adding rows consisting
of the $g$-vectors of $A(B^\dagger)$. (Note this unexpected fact:
the rows we add to $B$ are $g$-vectors, not for the cluster algebra $A(B)$,
but for the cluster algebra $A(B^\dagger)$. 
If 
$B$ is skew-symmetric, then $B^\dagger=-B$, and the $g$-vectors for $B^\dagger$
are the negatives of the $g$-vectors for $B$.

We are treating the above as the definition of universal
coefficients, but it is more properly a theorem of Nathan Reading
\cite[Corollary 8.15]{Re} that this produces what is defined as the cluster algebra with
universal coefficients in \cite{CA4}. (The definition from \cite{CA4} will not
be relevant for us.)

It is an important fact about universal coefficients that
$\mu_k(B^\univ)=\mu_k(B)^\univ$; that is to say, the coefficient rows of
$\mu_k(B^\univ)$ are exactly the rows we would obtain as above from the matrix
$\mu_k(B)$, i.e., the $g$-vectors of the cluster algebra $A(\mu_k(B)^\dagger)$
calculated with respect to the initial seed $\mu_k(B)^\dagger$.

In fact, we can be a little more precise. Let $z^\vee_1,\dots,z^\vee_v$ be the
cluster variables of $A(B^\dagger)$ in some order, and let $g^\vee_i$ be the
$g$-vector of $z^\vee_i$. 
It follows from the mutation rule for $g$-vectors under change of starting seed that when we apply $\mu_k$ to
$B^\univ$, the row $g^\vee_i$ is changed into the
$g$-vector of the same variable $z^\vee_i$, but considered with respect to the
mutated seed \cite[Conjecture 7.12]{CA4}, established in \cite[Proposition 11.3]{De}.
This says in particular
that there is a well-defined bijection between coefficient
rows of $B^\univ$ and cluster variables of $A(B^\dagger)$, which does not depend
on the choice of seed.

\subsection{An example}\label{ex}
For example, let $B$ be $\left[ \begin{array}{rr} 0 &-1\\2&0 \end{array} \right]$ and let $B^\dagger$ be
$\left [ \begin{array}{rr} 0 &2\\-1&0 \end{array} \right]$.
The matrix $(B^\dagger)^\prin$ is $$\left [ \begin{array}{rr} 0 &2\\-1&0 \\ \hline
    1&0\\0&1 \end{array} \right].$$

If we take our initial cluster variables for $A((B^\prin)^\vee)$ to be $x^\dagger_1$ and
$x^\dagger_2$, and our
frozen variables to be $y^\dagger_1$ and $y^\dagger_2$,
the other cluster variables are
\begin{eqnarray*}x_3^\dagger&=&\frac{x^\dagger_2+y^\dagger_1}{x^\dagger_1}\\
  x_4^\dagger&=&\frac{(x_2^\dagger)^2+2y^\dagger_1x^\dagger_2+(y_1^\dagger)^2+(y_1^\dagger)^2y^\dagger_2(x_1^\dagger)^2}{(x_1^\dagger)^2x^\dagger_2}\\
x_5^\dagger&=&\frac{x^\dagger_2+y^\dagger_1+y^\dagger_1y^\dagger_2(x^\dagger_1)^2}{x^\dagger_1x^\dagger_2}\\ x_6^\dagger&=&\frac{1+y^\dagger_2(x^\dagger_1)^2}{x^\dagger_2}.\end{eqnarray*}

The degrees of $x^\dagger_1$, $x^\dagger_2$, $y^\dagger_1$, $y^\dagger_2$ are
respectively
$(1,0)$, $(0,1)$, $(0,1)$, $(-2,0)$. The degrees of
$x_3^\dagger$, $x_4^\dagger$, $x_5^\dagger$, and $x_6^\dagger$ are then easily seen to be
$(-1,1)$, $(-2,1)$,$(-1,0)$, and $(0,-1)$. 

Therefore,
$$B^\univ=\left [ \begin{array}{rr} 0 &-1\\2&0 \\ \hline
    1&0\\0&1\\-1&1\\-2&1\\-1&0\\0&-1 \end{array} \right]$$

If we apply $\mu_1$ to $B^\univ$, we get 
$$\mu_1(B^\univ)=\left [ \begin{array}{rr} 0 &1\\-2&0 \\ \hline
    -1&0\\0&1\\1&0\\2&-1\\1&-1\\0&-1 \end{array} \right]$$
It is straightforward to check that these rows are indeed the $g$-vectors
of $\mu_1(B)^\dagger$, but this requires calculating them, which we shall skip.
It is easier
to check the corresponding assertion for $\mu_2\mu_1(B^\univ)$, because
$\mu_2\mu_1(B)=B$. Checking it, we find that

$$\mu_2\mu_1(B^\univ)=\left [ \begin{array}{rr} 0 &-1\\2&0 \\ \hline
    -1&0\\0&-1\\1&0\\0&1\\-1&1\\-2&1 \end{array} \right]$$
We observe that the rows are indeed the $g$-vectors of $B^\dagger$, which we had
already calculated. 

\subsection{Associahedra}\label{assoc}

Fixing, as before, a skew-symmetrizable matrix $B$, 
there is a fan in $\mathbb R^n$ called the $g$-vector fan, and denoted $\Sigma(B)$, whose rays are given by the
$g$-vectors of $A(B)$. A collection of $g$-vectors generates a cone of this fan
if the corresponding cluster variables live in a cluster together.

If the number of clusters of the cluster algebra is finite, then this fan
is complete, in the sense that its cones cover all of $\mathbb R^n$.

Recall that there is another way to produce a complete fan in $\mathbb R^n$: given a
full-dimensional
polytope in $\mathbb R^n$, take the rays to be generated by the normal vectors pointing
away from the facets, and define a collection of rays to generate a cone of
the fan
provided that they correspond to the facets meeting at some face of the polytope. The resulting fan
is called the outer normal fan of the polytope.

Given a complete fan, it is natural to ask
whether or not it is the outer normal fan
of a polytope. The question whether the $g$-vector fan could be realized in
this way goes back to the early days of the theory of cluster algebras
(and in fact, versions of this question were asked and answered before $g$-vectors were
defined). It is now standard to refer to a polytope realizing the
$g$-vector fan
in this way as a (generalized) associahedron.

There were a number of constructions given of realizations of
generalized associahedra. (See \cite{HPS} and references given there.) 
A strikingly simple realization of the associahedron
was given for a particular $B$-matrix of type $A$ in \cite{ABHY};
it was then shown that
the same construction could be applied for any simply laced finite type
cluster algebra and acyclic starting seed in \cite{BDMTY}. The restriction
on the starting seed was removed in \cite{PPPP,B}; \cite{PPPP} also
treats the skew-symmetrizable case.

We will recall the construction from \cite{PPPP} now. Note that they have
a non-standard convention for the entries of the $B$-matrix: their $B$-matrices are transposed compared to usual $B$-matrices. (Compare for example
their cluster exchange formula in section 2.2.1 with (4.2) of \cite{CA1}.)
This means that, when reproducing their formulas, we must introduce a
transposition compared to what is written in \cite{PPPP}.

Let $B$ be a skew-symmetrizable $n\times n$ matrix. Consider the cluster algebra
$A(B^\vee)$, with cluster variables $z^\vee_1,\dots,z^\vee_v$ and corresponding $g$-vectors
$g^\vee_1,\dots,g^\vee_v$. Suppose that we have numbered the cluster variables so that
$z^\vee_1,\dots,z^\vee_n$ are the initial cluster variables corresponding to the
starting matrix $B^\vee$.

We will record clusters as subsets $C$ of $\{1,\dots,v\}$ of size $n$. For
$C$ a cluster, we will write $B(C)$ for the corresponding $B$-matrix,
where we index the rows and columns of $B(C)$ by the entries of $C$, not
by $\{1,\dots,n\}$. 

Two cluster variables $z^\vee_i,z^\vee_j$ are called \emph{exchangeable} if there is a
cluster containing $z^\vee_i$, such that mutation at $z^\vee_i$ yields $z^\vee_j$.
A mutation from $z^\vee_i$ to $z^\vee_j$ is called \emph{positive mesh} if, when
we fix a cluster $C$ containing $z^\vee_i$, such that the result of mutating it is $z^\vee_j$,
then the entries $B(C)_{ki}$ for $k\in C$ 
are all non-negative. The property of being a positive mesh mutation does
not depend on which cluster containing $z^\vee_i$ is chosen among those such that the mutation of $z^\vee_i$ is $z^\vee_j$.

For each cluster
variable $z^\vee_i$, there is a unique $z^\vee_j$ which is the positive mesh mutation
of $z^\vee_i$. (For those familiar with the cluster category, we explain that
if $z^\vee_i$ corresponds to the indecomposable object $M$ in the cluster category,
then $z^\vee_j$
corresponds to the Auslander--Reiten translation of $M$. See \cite[Section 3]{PPPP} for more details.)

For any $i$ with $n<i\leq v$, let $z^\vee_j$ be the positive mesh mutation of $z^\vee_i$.
Then we have a relation among $g$-vectors:
$$g^\vee_i+g^\vee_j -\sum_{k\in C} B(C)_{ki}g^\vee_k = 0$$

(Note that if $1\leq i \leq n$, then there is still a well-defined $z^\vee_j$ such that $z^\vee_j$ is the positive mesh mutation of $z^\vee_i$, but the above relation does not hold. In fact, in this case, we simply have $g^\vee_i+g^\vee_j=0$.)

Further, these linear dependencies among the $g$-vectors provide a basis for
the space of linear dependencies among all the $g$-vectors. 

In our running example, for $3\leq i\leq 6$, the positive mesh mutation of
$x_i^\vee$ is to $x_{i-2}^\vee$.
The linear dependencies provided by the lemma are
$g^\vee_1-g^\vee_2+g^\vee_3=0$, $g^\vee_2-2g^\vee_3+g^\vee_4=0$, $g^\vee_3-g^\vee_4+g^\vee_5=0$, and $g^\vee_4-2g^\vee_5+g^\vee_6=0$.

Consider $\mathbb R^v$, with coordinates $w_1,\dots,w_v$.
For each $i$ with $n<i\leq v$, fix a positive
constant $c_i$.

For each $i$ satisfying $n<i\leq v$, let $z_j$ be the positive mesh mutation
of $z_i$, let $C$ be a cluster of $A(B^\dagger)$
with $i\in C$, such that the mutation of
$z_i$ is $z_j$, and let $B(C)$ be its corresponding $B$-matrix. Consider the following affine hyperplane in $\mathbb R^v$:

$$ w_j + w_{i} - \sum_{k\in C\setminus i} B(C)_{ki} w_k =c_i$$

Define $\mathbb E_c$ to be the intersection of all these affine hyperplanes.
Define $\mathbb U_c$ to be the intersection of $\mathbb E_c$ with the
non-negative orthant in $\mathbb R^v$.
Finally, define a projection $\pi:\mathbb R^v \rightarrow \mathbb R^n$,
by defining
$\pi$ to be the projection onto the coordinates $\tau(i)$, for $i$ initial.

Then \cite{PPPP} shows that
$\mathbb A_c=\pi(\mathbb U_c)$ is an associahedron, and all associahedra (up to
translation) can be obtained by suitable choice of $c$ (with each $c_i$ positive, as we have already assumed). In fact, as we
shall see, it is the polytope $\mathbb U_c$ which will naturally arise for
us in the context of symplectic reduction. ($\mathbb U_c$ is a polytope of the same dimension as $\mathbb A_c$, but embedded in a higher-dimensional space.)

In our running example, 
the equations cutting out $\mathbb U_c$ are as follows:

\begin{eqnarray*}
  w_1-w_2+w_3&=& c_3 \\
  w_2-2w_3+w_4&=& c_4\\
  w_3-w_4+w_5&=& c_5 \\
  w_4-2w_5+w_6&=&c_6
\end{eqnarray*}

\subsection{Degenerate associahedra as Newton polytopes}
Associahedra are important, among other things,
because degenerate associahedra give the Newton polytopes of $F$-polynomials.
We briefly give some further details here. This section is primarily motivational.

Let $B$ be a
skew-symmetrizable $n\times n$ matrix, such that the cluster algebra
$A(B)$ has finitely many cluster variables $z_1,\dots,z_v$, where we fix the
numbering so that $z_i=x_i$ for $1\leq i \leq n$ are the initial cluster variables.

The $F$-polynomials are defined as follows \cite{CA4}. Let $\widehat z_i$ be the cluster variable corresponding to $z_i$ in $A(B^\prin)$. These cluster
variables can be expressed as Laurent polynomials in $x_1,\dots,x_n,y_1,\dots,y_n$, with only the variables $x_j$ in the denominator. Then
$F_i$ is obtained from $\widehat z_i$ by substituting 1 for each $x_j$. The resulting $F_i$ is a polynomial in the variables $y_1,\dots,y_n$. For $1\leq i\leq n$, since $\widehat z_i=x_i$, we have $F_i=1$.

The Newton polytope of an $F$-polynomial, $\New(F_i)$, is by definition the convex hull of
the exponent vectors of the monomials appearing with non-zero coefficients in
the $F$-polynomial. Newton polytopes of $F$-polynomials have been a subject
of considerable recent interest, see \cite{BDMTY,PPPP,F1,F2,LLS,B,AHL}.

The Newton polytope of $F_i$ can be obtained as follows. Recall that
$\mathbb A_c$ depends on a $(v-n)$-tuple $c=(c_j)_{m<j\leq v}$. For
$n<i\leq v$, let $e_i$ denote the $(v-n)$-tuple which has a 1 in position $i$
and zeros elsewhere. Then for $n<i\leq v$, $\New(F_i)=\mathbb A_{e_i}$.

This was established under the assumption that $B$ is skew-symmetric and acyclic in \cite{BDMTY}. It was shown in \cite{B} that the acyclic assumption could be dropped. It was separately shown in \cite{AHL} via a folding argument that the skew-symmetric assumption could be dropped. By first unfolding, as in \cite{AHL},
and then applying the result of \cite{B}, the general result can be proved.

The polytope $\mathbb A_{e_i}$ is not itself an associahedron, in that its
outer normal fan is a coarsening of the $g$-vector fan. However, because the
coefficients of the $F$-polynomials are known to be 1 at the lattice points corresponding to their vertices, it is easy to see that multiplication of $F$-polynomials corresponds to Minkowski sum of Newton polytopes. Thus, any (possibly degenerate) associahedron $\mathbb A_c=\sum_i c_i\mathbb A_{e_i}$ with $c_i\in \mathbb Z_{\geq 0}$ is the Newton polytope of a suitable product of cluster variables.

\subsection{Symplectic reduction}
Let $X$ be a variety equipped with a real symplectic form $\omega$. Suppose
that a real torus $R=(S^1)^p$ has a Hamiltonian action on $X$ with respect to $\omega$.
There is a construction known as \emph{symplectic reduction} with
respect to a torus with a Hamiltonian action. Suppose that $R=R_1\times R_2$.
We will be interested in carrying out symplectic reduction with respect
to the subtorus $R_2$. 

Let $\ttt$ be the Lie algebra of $R$. (Since $R$ is commutative, $\ttt$ is
just $\mathbb R^p$ equipped with the zero Lie bracket.) Similarly, write
$\ttt_1,\ttt_2$ for the Lie algebras of $R_1$ and $R_2$.
There is a moment map for the $R$ action, which we call $\mu_{R}$, sending
$X$ to $\ttt^*$ (i.e., linear functions on $\ttt$). By
restriction, this gives us a map $\mu_{R_2}:X\rightarrow \ttt_2^*$,
where $\mu_{R_2}=i^*\circ \mu_{R}$ with $i^*$ being
the map induced from the inclusion of $\ttt_2$ into $\ttt$. 

Fix $\hat c\in \ttt_2^*$. Assume that $R_2$ acts freely on $\mu_{R_2}^{-1}(\hat c)$.
Then the quotient $X_\red=\mu_{R_2}^{-1}(\hat c)/R_2$ is known as the symplectic reduction
of $X$ with respect to $R_2$. The Marsden--Meyer--Weinstein theorem (see,
for example, \cite[Chapter 23]{LSG}) says that the symplectic reduction inherits a symplectic form $\omega_\red$.

$R/R_2\simeq R_1$ has a Hamiltonian action on $X_\red$, and by
\cite[Remark, p.~146]{LSG}, we have the following proposition:

\begin{proposition} \label{intersect} The moment map image of $X_\red$ with respect to the
action of $R_1$ is given by
$(i^*)^{-1}(\hat c) \cap \mu_R(X)$, where we identify $(i^*)^{-1}(\hat c)$ with $\ttt^*_1$.
\end{proposition}

To rephrase this proposition in words, the moment map image of $X_\red$ with respect to the action of $R_1$ is given by taking an affine slice of the moment map image of $X$ with respect to the action of $R$. The alert reader will note a resemblance to the construction of associahedra from \cite{PPPP}.

\subsection{Cluster varieties} \label{varieties}

A cluster algebra $A(\widetilde B)$ defined an integral $(n+m)\times n$ matrix $\widetilde B$ of full rank with a skew-symmetrizable principal part can be identified
with a ring of regular functions on a nonsingular rational Poisson variety $\mathcal{M}(\widetilde B)$, 
first defined in \cite{GSV1} and traditionally called $\A$-variety.   Coordinate charts of this variety correspond to clusters in $A(\widetilde B)$. 
Consider a cluster $s=\{ A_1^{(s)},\dots, A_{n+m}^{(s)}\}$. (This is just another notation for the cluster variables with the last $m$ variables being frozen.)
The corresponding coordinate chart is $U^{(s)}=\Spec({A_1^{(s)}}^{\pm 1},\dots, {A_n^{(s)}}^{\pm 1}, A_{n+1}^{(s)},\dots, A_{n+m}^{(s)})$. If the cluster $s'$ is
obtained from cluster $s$ by mutation at position $i$, then the same mutation defines the transition map between $U^{(s)}$ and $U^{(s')}$. That is to say, if two clusters are related by mutation in position $i$, then
$A^{(s')}_j=A^{(s)}_j$ for $j\ne i$, and $A^{(s')}_i$ is obtained from the
mutation formula.

$\mathcal{M}(\widetilde B)$ is equipped with an action of an $m$-dimensional torus $(\mathbb{C}^*)^m$ which can be described as follows.
For an arbitrary  cluster
$s=\{ A_1^{(s)},\dots, A_{n+m}^{(s)}\}$ we define a {\it local toric action\/} of rank $m$ as the map 
$T_K^{\d}:\FF_\C\to
\FF_\C$ given on the generators of $\FF_\C=\C(A_1^{(s)},\dots, A_{n+m}^{(s)})$ by the formula 
\[
T_K^{\d}(A_i^{(s)})=A_i^{(s)}\prod_{\alpha=1}^m d_\alpha^{k_{i\alpha}},\quad i\in [n+m],\qquad
\d=(d_1,\dots,d_r)\in (\C^*)^m,
\label{toricact}
\]
where $K=K^{(s)}=(k_{i\alpha})$ is an integer $m\times (n+m)$ {\it weight matrix\/} of full rank, and extended naturally to the whole $\FF_\C$. 

Let $s'$ be another cluster, then the corresponding local toric action defined by the weight matrix $W'$
is {\it compatible\/} with the local toric action if the following diagram is commutative for
any fixed $\d\in (\C^*)^r$:
$$
\begin{CD}
\FF_\C=\C(s) @>>> \FF_\C=\C(s')\\
@V T_K^{\d} VV @VV T_{K'}^{\d}V\\
\FF_\C=\C(s) @>>> \FF_\C=\C(s')
\end{CD}
$$
Here the horizontal arrows are induced by mutations $A_i^{(s)}\mapsto A_i^{(s')}$. If local toric actions at all clusters are compatible, they define a {\it global toric action\/} $T^{\d}$ on $\FF_\C$ called the extension of the local toric action above. Lemma~2.3 in \cite{GSV1} claims that such extension 
to a unique global action of $(\C^*)^m$  exists if and only if $K \widetilde B^{(s)}= 0$. Thus the left kernel of the initial exchange matrix $\widetilde B$ defines a toric action of $(\mathbb{C}^*)^m$ on $\mathcal{M}(\widetilde B)$.

\section{Proof of the main result}

\subsection{Symplectic reduction of a patch of $\mathcal M(B^\univ)$}
Fix a skew-symmetrizable $n\times n$ matrix $B$, which corresponds to
a finite type cluster algebra with $v$ cluster variables, and let $B^\univ$ be the same $n\times n$
matrix, with $v$ rows added, providing universal coefficients.

The patch of $\mathcal M(B^\univ)$ corresponding to this cluster looks like
$P=(\mathbb C^*)^n \times \mathbb C^v$. The torus $T=(\mathbb C^*)^{n+v}$ acts
naturally on it.

We put a real symplectic form on $P$. Note that this is not the holomorphic
log-canonical pre-symplectic form compatible with the cluster structure,
as defined, for example, in \cite{GSV2}.
For $1\leq j \leq n+v$, let $z_j=r_j(\cos \theta_j+i\sin \theta_j)$ be coordinates on $P$, with $z_1,\dots,z_n$ corresponding to the initial coordinates in $\mathbb C^*$. 
Let $$\omega=\sum_{j=1}^n \frac{dr_j\wedge d\theta_j}{r_j} + \sum_{j=n+1}^{n+v} -r_jdr_j \wedge d\theta_j.$$

Consider the vector field $X_j=\partial/\partial \theta_j$. For
$1\leq j\leq n$,
$$i_{X_j}\omega=-dr_j/r_j=-d\log(r_j).$$
For $n+1\leq j \leq n+v$,
$$i_{X_j}\omega=r_jdr_j=\frac12 dr_j^2.$$

This shows that the action of the torus is Hamiltonian, and
the moment map is given by $\mu(z_1,\dots,z_{n+v})=(-\log(r_1),\dots,-\log(r_n),\frac12 r_{n+1}^2, \dots,\frac 12 r_{n+v}^2).$
It follows that the moment map image is $\mathbb R^n \times \mathbb R_{\geq 0}^v$.

Using the notations of Section \ref{varieties}, denote by $K$ an $v\times (n+v)$ integer matrix whose rows form a basis the left kernel of $\widetilde B$. 
We will pick a specific form for $K$ as follows: its first $n$ rows correspond to expressing the first $n$ rows of $\widetilde B$ as linear combinations of the final $v$ rows and the remaining $(v-n)$ columns reflect 
linear dependencies among the final $v$ rows of $\widetilde B$.


Now, we  consider symplectic reduction with respect to the real part  $T_K(\mathbb{R})$ of the torus $T_K$. We  view  $T_K(\mathbb{R})$ as a subtorus of  $T(\mathbb R)=(S^1)^{n+v}$ whose action on the patch $P$ is described above. To this end, we consider a restriction $\mu_K$ of the moment map $\mu_T$ defined as in Section 2.6. Pick $\hat c\in \mathbb R^n \times \mathbb{R}_{>0}^{v-n}$. Write $c$ for the final $(v-n)$ entries of $\hat c$. The fact that the entries of $c$ are positive guarantees that $T_K(\mathbb{R})$ acts freely on $\mu_K^{-1}(\hat c)$. We are interested in the moment map image of the  quotient $\mu_{K}^{-1}(\hat c)/T_K(\mathbb{R})$ with respect to the action of $T/T_K(\mathbb{R})$.
The answer, given by Proposition \ref{intersect}, can be made more explicit due to our choice of $K$. Indeed, 
in carrying out the symplectic reduction, we intersect $\mathbb R^n \times \mathbb R_{\geq 0}^v$ with $(i^*)^{-1}(\hat c)$. 
Comparing our description of the last $k$ columns of $K$ with the construction outlined Section 2.3, we recognize that in computing this intersection we use precisely the 
equations defining the affine slice that cuts out the associahedron. The way the first $n$ columns of $K$ are defined
 allows us to solve for the coordinates in $\mathbb R^n$ in terms of the
other data, and since these coordinates do not have any constraints on them,
we have ended up with the ABHY associahedron associated to $c$. We summarize this in the following
theorem:

\begin{thm} Consider the symplectic reduction of
  $\mathcal M(B^\univ)$ by the action of $T_K(\mathbb R)$ at
  $\hat c\in \mathbb R^n \times \mathbb R^{v}_{>0}$.
  The moment map image of the symplectic reduction
  yields the associahedron $\mathbb U_c$ associated to $c$
  by \cite{PPPP}.
  \end{thm}

\begin{remark}
  {\rm The \emph{cluster configuration space} studied by Arkani-Hamed, He, and
    Lam \cite{AHL} can be realized as a dense set inside the
    quotient of $\mathcal M(B^\univ)$ by $T_K$.
    One easily verifies that symplectic reduction of
    $\mathcal M(B^\univ)$ by $T_K(\mathbb R)$ amounts to the same thing as taking the quotient by $T_K$. We view this as providing some motivation for the naturalness of the symplectic quotient considered here.}
\end{remark}

\begin{remark} \rm The fact that the associahedron is defined as the moment
  map image with respect to a torus which acts only on $P$, not on
  $\mathcal M(B^\univ)$, explains the fact that different choices of initial seed yield different associahedra. \end{remark}

\subsection{An example}
Let us consider the rank two example which we have already discussed.

Here,
$$B^\univ=\left [ \begin{array}{rr} 0 &-1\\2&0 \\ \hline
    1&0\\0&1\\-1&1\\-2&1\\-1&0\\0&-1 \end{array} \right]$$

The cluster algebra $A(B^\univ)$ is naturally
$\mathbb Z^v$-graded; the grading lattice is more properly the
left kernel of $B^\univ$.

A basis for the left kernel is given by the rows of
the following matrix:

$$K=\left[\begin{array}{rr|rrrrrr}
        1&0&0&1&0&0&0&0\\
0&1&-2&0&0&0&0&0\\
    0&0&1&-1&1&0&0&0\\
    0&0&0&1&-2&1&0&0\\
    0&0&0&0&1&-1&1&0\\
  0&0&0&0&0&1&-2&1\end{array}\right]$$

The top two rows express the first two rows of $B^\univ$ as a linear combination of the
$g$-vectors, while 
the bottom four rows are the linear dependencies among the $g$-vectors.

We pick a 6-tuple $\hat c =(c_1,\dots,c_6)$ in $\mathbb R^2 \times \mathbb R_{>0}^4$, and we intersect
$\mathbb R^2 \times \mathbb R^6_{\geq 0}$ with the following conditions, where we write $(u_1,u_2)$ for the coordinates in $\mathbb R^2$ and
$(w_1,\dots,w_6)$ for the coordinates in $\mathbb R^6_{\geq 0}$. 

\begin{eqnarray*}   u_1+w_2 &=& c_1\\
  u_2 - 2w_1&=&c_2\\
  w_1-w_2+w_3&=&c_3\\
  w_2-2w_3+w_4&=&c_4\\
  w_3-w_4+w_5&=&c_5\\
w_4-2w_5+w_6&=&c_6\end{eqnarray*}

  The first two equations fix $u_1,u_2$ in terms of the other data, and since
  there are no constraints on $u_1,u_2$, these coordinates can be forgotten. The
  remaining equations are exactly the equations that we saw at the end of
  Section \ref{assoc} for the associahedron associated to $c=(c_3,c_4,c_5,c_6)$.

  \section*{Acknowledgements} H.T. thanks Ben Webster for an excellent
symplectic geometry course, for which this was his final homework project
(on which it turned out he needed assistance). The authors also thank
Jacques Hurtubise, Allen Knutson, and Milen Yakimov for their valuable
suggestions.
M.G. was partially supported by NSF grant DMS \#2100785. H.T. was partially
supported by NSERC Discovery Grant RGPIN-2022-03960 and the
Canada Research Chairs program, grant number CRC-2021-00120.
The authors express their gratitude to the Mathematisches
Forschingsinstitut Oberwolfach, where this paper was finished at Workshop
2403.

\end{document}